\input amstex
\documentstyle {amsppt}

\pagewidth{32pc} 
\pageheight{45pc} 
\mag=1200
\baselineskip=15 pt

\hfuzz=5pt 
\topmatter
\NoRunningHeads 
\title Appendix to Discrete Localities II: $p$-local compact groups as localities  
\endtitle
\author Andrew Chermak and Alex Gonzalez
\endauthor 

\address Mathematics Department, Kansas State University, Manhattan, KS 66506 (USA) 
\endaddress 

\email chermak\@math.ksu.edu 
\endemail 

\address Institut Arraona, 
43 Praga Street, 
08207 Sabadell (Spain)
\endaddress 

\email agondem\@gmail.com
\endemail 

\date
March 2022 
\enddate 

\endtopmatter

\redefine\norm{\trianglelefteq}

\redefine\bar{\overline}

\redefine\maps{\mapsto}
\redefine\i{^{-1}}

\redefine\l{\lambda}
\redefine\s{\sigma}
\redefine\a{\alpha}
\redefine\b{\beta}
\redefine\d{\delta}
\redefine\g{\gamma}
\redefine\e{\epsilon} 

\redefine\t{\tau}
\redefine\r{\rho}

\redefine\G{\Gamma}

\redefine\<{\langle}
\redefine\>{\rangle}

\redefine\ca{\Cal}

\redefine\D{\Delta}

\redefine\sub{\subseteq}

\redefine\1{\bold 1} 

\redefine\up{\uparrow} 
\redefine\down{\downarrow}

\redefine\bull{\bullet}

\document 
   
Our aim in this appendix is to show that the $p$-local compact groups, introduced in [BLO2] and 
further developed in [BLO3], [JLL] and elsewhere, may be viewed as proper localities of a certain kind 
- to be called compact localities. As a corollary to a result of Ran Levi and Assaf Libman [LL] (and 
with an improvement due to Remi Molinier [M]), we will show that if $\ca F$ is the fusion system of a 
compact locality $(\ca L,\D,S)$ then, up to an isomorphism of partial groups which restricts to an  
automorphism of $S$, $\ca L$ is the unique compact locality on $\ca F$ having $\D$ as its set of objects. 
This result, and the formalism of compact localities, provide a bridge between the theory being developed 
in this series and a more well-established theory having a much more homotopy-theoretic flavor. 

We shall be closely following the arguments in the appendix 
to [Ch1], where an equivalence was established between the ``$p$-local finite groups" 
introduced in [BLO1], and proper, finite, centric localities. 

Recall from the Appendix A in Part I (or from any of the above references) that a {\it $p$-torus} $T$ is by 
definition the direct product of a finite number of copies of the Pr\" ufer group $\Bbb Z/{(p^\infty)}$. 
A group $G$ is {\it virtually $p$-toral} if there exists a $p$-torus of finite index in $G$. 
A virtually $p$-toral $p$-group is a {\it discrete $p$-toral group}. 

\definition {Definition A.1} The discrete locality $(\ca L,\D,S)$ on $\ca F$ is {\it compact} if:  
\roster 

\item "{(1)}" $\ca L$ is proper, and 

\item "{(2)}" $N_{\ca L}(P)$ is viritually $p$-toral for each subgroup $P\in\D$. 

\endroster 
\enddefinition 

We will show via A.7, A.22, and A.23 below that there is an equivalence between the notions of  
compact locality and $p$-local compact group. We then obtain an existence/uniqueness theorem for 
compact localities, as a corollary to [LL]. 

\vskip .1in 
The set of elements $x$ in a $p$-torus $T$ such that 
$x^p=\1$ is an elementary abelian $p$-group of finite order $p^k$ where $k$ is equal to 
the number of factors in any decomposition of $T$ as a direct product of Pr\" ufer groups. We 
refer to $k$ as the {\it rank} of $T$, and write $rk(T)=k$.  
The identity group is the $p$-torus of rank $0$. 

\proclaim {Lemma A.2} Let $G$ be a virtually $p$-toral group and let $P$ be a discrete $p$-toral group. 
\roster 

\item "{(a)}" There is a unique $p$-torus $T$ such that $T$ has finite index in $G$; and then $T$ contains 
every $p$-toral subgroup of $G$. 

\item "{(b)}" All subgroups and homomorphic images of $G$ are virtually $p$-toral,  
and all subgroups and homomorphic images of $P$ are discrete $p$-toral. 

\item "{(c)}" If $X$ and $Y$ are subgroups of $P$ with $X<Y$ ($X$ is a proper subgroup of $Y$) then 
$X< N_Y(X)$.  

\endroster 
\endproclaim 

\demo {Proof} Let $X$ be a $p$-torus contained in $G$. Then $T\cap X$ has finite index in $X$, and then 
$T\cap X=X$ since $X$ is $p$-divisible. This proves (a). Now let $H$ be a subgroup of $G$. As   
$T\norm G$ by (a), there is an isomorphism $HT/T\cong H/(H\cap T)$, and thus $H$ is discrete 
$p$-toral. Let $N\norm G$. Then $TN/N$ has finite index in $G/N$, and $TN/N\cong T/(N\cap T)$ where 
$T/(N\cap T)$ by $p$-divisibility. This proves (b). Point (c) is given by Lemma A.3(a) in 
the Appendix to Part I. 
\qed 
\enddemo 

We refer to the unique $p$-torus of finite index in the virtually $p$-toral group $G$ as the 
maximal torus of $G$. If $P$ is a discrete $p$-toral group with maximal torus $T$ then (following  
[BLO2]) the {\it order} of $P$ is defined to be the pair 
$$ 
|P|=(rk(P),|P/T|). 
$$ 
If $Q$ is a discrete $p$-toral group with maximal torus $U$ then we write  
$$ 
|P|<|Q| 
$$ 
if $(rk(P),|P/T|)<(rk(Q),|Q/U|)$ lexicographically (i.e. $rk(P)<rk(Q)$, or $rk(P)=rk(Q)$ and $|P/T|<|Q/T|$).

\proclaim {Corollary A.3} Let $(\ca L,\D,S)$ be a compact locality on $\ca F$, and let $(\ca L',\D',S)$ be 
an expansion or a restriction of $\ca L$ to a proper locality. Then $\ca L'$ is a compact locality 
on $\ca F$. 
\endproclaim 

\demo {Proof} By 7.2(a) $\ca L'$ is a proper locality on $\ca F$. Thus we need  
only verify that the condition (2) in the preceding definition holds with 
$(\ca L',\D',S)$ in the place of $(\ca L,\D,S)$. Here (2) is obvious if $\D'\sub\D$. By Zorn's Lemma 
it suffices to consider the case where $(\ca L',\D',S)$ is an elementary expansion of $(\ca L,\D,S)$. 
That is, with $\ca F=\ca F_S(\ca L)$, we may assume that $\D'=\D\cup R^{\ca F}$ for some $R\leq S$,  
and we may then take $R$ to be fully normalized with 
respect to the stratification on $\ca F$ induced from $\ca L$. Then $N_{\ca L'}(R)$ is a subgroup of 
$N_{\ca L}(Q)$ for some $Q\in\D$ by 7.1, and then A.2(b) shows that $N_{\ca L'}(R)$ is virtually 
$p$-toral, as required. 
\qed 
\enddemo

\proclaim {Lemma A.4} Let $P$ and $Q$ be discrete $p$-toral groups. 
\roster 

\item "{(a)}" If $P\cong Q$ then $|P|=|Q|$. 

\item "{(b)}" If $P\leq Q$ then $|P|\leq |Q|$, with equality if and only if $P=Q$. 

\item "{(c)}" If $P\leq Q$ and $Q$ is a $p$-torus, then either $P=Q$ or $rk(P)<rk(Q)$. 

\item "{(d)}" Assume $P\leq Q$, set $P_0=P$, and recursively define $P_k$ for $k>0$ by 
$P_k=N_Q(P_{k-1})$. Set $B=\bigcup\{P_k\}_{k\geq 0}$. Then either $B=Q$ or $rk(B)<rk(P)$. 

\endroster 
\endproclaim  

\demo {Proof} Points (a) and (b) are straightforward. Now let $P\leq Q$ and assume that $Q$ is a 
$p$-torus. For any abelian $p$-group $A$ and any $m>0$ let $A_m$  
be the subgroup of $A$ generated by elements of order dividing $p^m$. Assuming that $P$ is 
a proper subgroup of $Q$, we find $Q_m\nleq P$ for $m$ sufficiently large, and then the rank of 
$P_m$ is less than the rank of $Q_m$. This yields $rk(P)<rk(Q)$, and thus (c) holds. 

In proving (d) let $V$ be the maximal torus of $B$, and let $V_0$ the 
maximal torus of $P_0$. Then $V_0\leq V$. Suppose that $rk(B)=rk(P)$. Then $V_0=V$, so $|B:P|$ is 
finite, and so there exists $k$ with $P_k=P_{k+1}=B$. Then $B=P$ by A.2(c), and thus (d) holds.  
\qed 
\enddemo 

Before we can begin working towards a correspondence between compact localities and $p$-local 
compact groups (and before we can state the definition of $p$-local compact group), we have first to sort 
out some potentially conflicting terminology 
concerning fusion systems. The difficulty here stems from our having already established some terminology in 
sections 2 and 8 of Part II, relating to stratified fusion systems. This terminology will be shown to 
agree with  that of the cited references on $p$-local compact groups 
once the context has been narrowed down to the fusion systems associated with linking 
systems, but until that point has been reached there is a very real problem of confusion.  
For example, if $(\ca F,\Omega,\star)$ is a stratified fusion system on a $p$-group $S$ 
then we have the notion (II.2.7) of a subgroup $P\leq S$ being fully normalized or fully centralized  
in $\ca F$, but there is a quite different definition of these terms in [BLO2]; and there is a 
similar difficulty concerning saturation of fusion systems. Our solution 
is to slightly alter the terminology from [BLO2], in the manner of the following definition.

\definition {Definition A.5} Let $\ca F$ be a fusion system over the discrete $p$-toral group $S$. 
\roster 

\item "{$\cdot$}" A subgroup $P\leq S$ is {\it fully order-centralized in $\ca F$} if
$|C_S(P)|\geq|C_S(Q)|$ for all $Q\in P^{\ca F}$. 

\item "{$\cdot$}" A subgroup $P\leq S$ is {\it fully order-normalized in $\ca F$} if 
$|N_S(P)|\geq|N_S(Q)|$ for all $Q\in P^{\ca F}$. 

\item "{$\cdot$}" $\ca F$ is {\it order-saturated} if the following three conditions hold. 

\item "{(I)}" For each $P\leq S$ the group $Out_{\ca F}(P)$ is finite. Moreover, if 
$P$ is fully order-normalized in $\ca F$, then $P$ is fully 
order-centralized in $\ca F$, and 
$$Out_S(P)\in Syl_p(Out_{\ca F}(P)).$$  

\item "{(II)}" If $P\leq S$ and $\phi\in Hom_{\ca F}(P,S)$ are such that $P\phi$ is fully order-centralized 
in $\ca F$, and if we set 
$$ 
N_{\phi}=\{g\in N_S(P)\mid \phi\i\circ c_g\circ\phi\in Aut_S(P\phi)\}, 
$$ 
then there exists $\bar\phi\in Hom_{\ca F}(N_\phi,S)$ such that $\phi=\bar\phi\mid_P$. 

\item "{(III)}" If $P_1<P_2<P_3<\cdots$ is an increasing sequence of subgroups of $S$, with 
$P_\infty=\bigcup\{P_n\}_{n=1}^\infty$, and if $\phi:P_\infty\to S$ is a homomorphism such that 
$\phi\mid_{P_n}\in Hom_{\ca F}(P_n,S)$ for all $n$, then $\phi\in Hom_{\ca F}(P_\infty,S)$. 

\endroster 
\enddefinition

By [BLO2, Lemma 1.6], if $\ca F$ is a fusion system over a discrete $p$-toral group $S$, then for any 
subgroup $P\leq S$ there are upper bounds for $|N_S(Q)|$ and for $|C_S(Q)|$ taken over $Q\in P^{\ca F}$. 
Thus $P$ has at least one fully order-normalized $\ca F$-conjugate and at least one fully order-centralized 
$\ca F$-conjugate.

\proclaim {Lemma A.6} Let $(\ca L,\D,S)$ be a compact locality on the order-saturated 
fusion system $\ca F$, and let $P\leq S$ be fully normalized in $\ca F$ with respect to the 
stratification induced from $\ca L$. 
\roster 

\item "{(a)}" $P$ is fully order-normalized in $\ca F$. 

\item "{(b)}" $Inn(P)=O_p(Aut_{\ca F}(P))$ if and only if $P=O_p(N_{\ca F}(P))$. 

\endroster 
\endproclaim 

\demo {Proof} As remarked above, there exists $Q\in P^{\ca F}$ such that 
$Q$ is fully order-normalized in $\ca F$. Set $U=N_S(P)$ and $V=N_S(Q)$. As $\ca F$ is 
order-saturated there then exists an $\ca F$-homomorphism $\phi:U\to V$ with $P\phi=Q$. 

Set $\Omega=\Omega_S(\ca L)$ and let $(\Omega,*)$ be the stratification on $\ca F$ induced from $\ca L$.  
For $X\leq S$ write $dim(X)$ for $dim_{\Omega}(X)$. As $P$ is fully normalized we have 
$dim(U)\geq dim(V)$, and then equality holds since $U\phi\leq V$. Let $\psi$ 
be an extension of $\phi$ to an $\ca F$-homomorphism $U^*\to V^*$. Then $U^*\psi=V^*$ since 
$dim(U^*)=dim(U)$ and $dim(V^*)=dim(V)$. Then also $U\phi=V$ since $U=N_{U^*}(U)$ and 
$V=N_{V^*}(V)$. Thus (a) holds. 

As $\ca F$ is order-saturated there exists a subgroup $\bar P$ of $N_S(P)$ such that $P\leq\bar P$ and 
$Aut_{\bar P}(P)=O_p(Aut_{\ca F}(P))$. Condition (II) in definition A.6 then yields 
$\bar P\norm N_{\ca F}(P)$. So, if $P=O_p(N_{\ca F}(P))$ then $P=\bar P$, and 
thus $Inn(P)=O_p(Aut_{\ca F}(P))$. Now set $R=O_p(N_{\ca F}(P))$ and notice that 
$Aut_R(P)\norm Aut_{\ca F}(P)$. Then $P=R$ if $P=\bar P$, and this completes the proof of (b). 
\qed 
\enddemo

From now on, when speaking of the fusion system $\ca F=\ca F_S(\ca L)$ of a compact locality, if we  
say that a subgroup $P\leq S$ is fully normalized in $\ca F$ then we mean that $P$ is fully normalized 
with respect to the stratification induced from $\ca L$; and similarly for ``fully centralized". 

\proclaim {Lemma A.7} Let $\ca F$ be a fusion system on the discrete $p$-toral group $S$. Assume: 
\roster 

\item "{(1)}" $\ca F$ is saturated (as defined in 8.2), and 

\item "{(2)}" $Out_{\ca F}(P)$ is finite for all $P\leq S$.

\endroster 
 Then $\ca F$ is order-saturated. Further: 
\roster 

\item "{(a)}" A subgroup $P\leq S$ is fully order-normalized in $\ca F$ if and only if $P$ is 
fully normalized in $\ca F$. 

\item "{(b)}" A subgroup $P\leq S$ is fully order-centralized in $\ca F$ if and only if $P$ is 
fully centralized in $\ca F$. 

\item "{(c)}" Let $T$ be the maximal torus of $S$.   
Let $R$ be any subgroup of $S$ such that $R$ is a $p$-torus, let $R'$ be an $\ca F$-conjugate of $R$, 
and let $\a:R\to R'$ be an $\ca F$-isomorphism. Then $R$ and $R'$ are subgroups of $T$, and $\a$ 
extends to an $\ca F$-automorphism of $T$. 

\endroster 
Moreover, the conditions (1) and (2) obtain if $\ca F$ is the fusion system of a compact locality. 
\endproclaim 

\demo {Proof} Let $P\leq S$ be fully 
order-normalized in $\ca F$, and let $Q\in P^{\ca F}$ be fully normalized. Then there exists an 
$\ca F$-homomorphism $\phi:N_S(P)\to N_S(Q)$ with $P\phi=Q$. 
As $|N_S(P)|\geq|N_S(Q)|$ it follows from A.4 that $N_S(P)\phi=N_S(Q)$. Thus $P$ is fully normalized in 
$\ca F$, and $Q$ is fully order-normalized in $\ca F$. This yields (a), and (b) follows in similar fashion. 

We next show that $\ca F$ satisfies the condition (I) in A.5. In view of the hypothesis (2) it need only 
consider the case where $P$ is fully order-normalized in $\ca F$. Then $P$ is fully normalized 
(by (a)), hence fully centralized (by I.2.9), and hence fully order-centralized (by (b)). 
By definition 8.2, $P$ is fully automized, so $Out_S(P)\in Syl_p(Out_{\ca F}(P))$ as required. The 
condition (II) in A.5 is satisfied since (by 8.2) $\ca F$ is receptive. 

Before verifying condition (III) in A.5, we shall need to prove (c). 
Let $R$ be any $p$-torus contained in $S$, let $T$ be the maximal torus of $S$, 
and let $\a:R\to R'$ be an $\ca F$-isomorphism. Then $R$ and $R'$ are subgroups of $T$ by A.2(a). 
Let $V\in R^{\ca F}$ be fully centralized in $\ca F$ and since $T\leq C_S(R)$ there exists an 
$\ca F$-automorphism $\b:T\to T$ with $R'\b=V$. Then $\a\circ\b$ is an $\ca F$-isomorphism 
$R\to V$, and $\a\circ\b$ extends to an $\ca F$-automorphism $\g$ of $T$ since $V$ is receptive. 
Then $\g\circ\b\i$ is an extension of $\a$ to $T$, and so (c) holds. 

We now turn to (III) in A.5. Let $\s=(P_i)_{i=1}^\infty$ be an increasing sequence of subgroups of 
$S$ and let $\t=(\phi_i)_{i=1}^\infty$ be a sequence of $\ca F$-homomorphisms, $\phi_i:P_i\to S$, such that 
$\phi_i=\phi_{i+1}\mid_{P_i}$ for all $i$. Let $P$ be the union of the groups $P_i$ and let 
$\phi:P\to S$ be the union of the mappings $\phi_i$. For each infinite subset $I$ of natural numbers, 
$\phi$ is then the union of the mappings $\phi_i$ for $i\in I$; so in order to prove that $\phi$ is an 
$\ca F$-homomorphism we are free to replace $\s$ and $\t$ by the sequences corresponding to $I$, and 
then to assume that $I=\Bbb N$. As $\ca F$ is stratified we have 
$(P_i)^*=P^*$ for $i$ sufficiently large, and so we may assume that  
$(P_i)^*=P^*$ for all $i$. Then each $\phi_i$ extends to an $\ca F$-homomorphism $P^*\to S$, which 
then restricts to an $\ca F$-homomorphism $\psi_i:P\to S$. 

Let $U$ be the maximal torus of $P$. Then $P=UP_n$ for $n$ sufficiently large, and so we may assume 
$P=UP_1$. As $Aut_{\ca F}(T)=Out_{\ca F}(T)$ is finite, (c) implies that there are only 
finitely many $\ca F$-conjugates of $U$, and indeed that there exists an $\ca F$-homomorphism 
$\xi:U\to S$ and an infinite set $I$ of natural numbers such that $\psi_i\mid_U=\xi$ for all $i\in I$. 
As $P=UP_1$ the set $\{\psi_i\}_{i\in I}$ then has a single element $\psi$, and $\psi$ is then the 
union of the homomorphisms $\phi_i$. Thus $\psi=\phi$, so $\phi$ is an $\ca F$-homomorphism, and 
(III) holds. 

Finally, assume that we are given a compact locality $(\ca L,\D,S)$. Then $\ca L$ is proper, so 
$\ca F_S(\ca L)$ is saturated by 8.3(c). Thus, it only remains to show that the condition (2) holds 
for $\ca F=\ca F_S(\ca L)$. By Theorem 7.2 there exists a proper expansion 
$(\ca L^+,\D^+,S)$ on $\ca F$, with $\ca F^c\sub\D^+$; and $\ca L^+$ is then compact by A.3. Thus, 
we may assume $\ca F^c\sub\D$. 

As $(\ca L,\D,S)$ is proper we may appeal to Theorem 7.2 to obtain a proper expansion 
$(\ca L^+,\D^+,S)$ on $\ca F$, with $\ca F^c\sub\D^+$. By construction, subgroups of $\ca L^+$ are 
conjugates of subgroups of $\ca L$, and thus (2) holds and $\ca L^+$ is compact. We may therefore 
replace $\ca L$ 
with $\ca L^+$ in the remainder of the proof. That is, we may assume that we have $\ca F^c\sub\D$. 
Let $P\leq S$ be fully normalized in $\ca F$. As $\ca F$ is inductive, $P$ is then fully centralized 
in $\ca F$.   Set $X=C_S(P)P$. Then $X$ is $\ca F$-centric, and so $X\in\D$. Set $H=N_{\ca L}(X)$. Every 
$\ca F$-automorphism of $P$ extends to an $\ca F$-automorphism of $X$ by receptivity, so 
$$ 
Aut_{\ca F}(P)=Aut_{N_H(P)}(P).   
$$   
Thus $Aut_{\ca F}(P)$ is a homomorphic image of the virtually $p$-toral group $N_H(P)$. In particular 
every element of $Aut_{\ca F}(P)$ is of finite order. By [BLO2, Lemma 1.5(b)] every torsion subgroup of 
$Out(P)$ is finite, so we conclude 
that $Out_{\ca F}(P)$ is finite. That is, (2) holds, and the proof is complete. 
\qed 
\enddemo

This completes the preliminaries concerning fusion systems. We next show how a compact locality 
$(\ca L,\D,S)$ on $\ca F$ gives rise to a ``transporter system" of a certain kind, and then that this 
transporter system is a $p$-local compact group provided that $\D$ is the set $\ca F^c$
of $\ca F$-centric subgroups of $S$. The definition of 
transporter system will be taken from [BLO3].

\vskip .1in 
Let $(\ca L,\D,S)$ be any locality. For each $(P,Q)\in\D\times\D$ let $N_{\ca L}(P,Q)$ be the set of 
all $g\in\ca L$ such that $P\leq S_g$ and $P^g\leq Q$. There is then a category $\ca T=\ca T_{\D}(\ca L)$ 
whose set of objects is $\D$, and whose morphisms $g:P\to Q$ are triples $(g,P,Q)$ such that 
$g\in N_{\ca L}(P,Q)$, with composition defined by 
$$ 
(g,P,Q)\circ (h,Q,R)=(gh,P,R). 
$$  
In practice, the role of the objects $P$ and $Q$ will always be clear from the context, and 
we may therefore identify $Mor_{\ca T}(P,Q)$ with $N_{\ca L}(P,Q)$. 

Let $\ca F$ be the fusion system $\ca F_S(\ca L)$, i.e. the fusion system on $S$ generated by the conjugation 
maps $c_g:P\to Q$ with $P,Q\in\D$ and with $g\in N_{\ca L}(P,Q)$. There is then a functor 
$$ 
\r:\ca T\to\ca F 
$$ 
such that $\r$ is the inclusion map $\D\to Sub(S)$ on objects, and such that 
$\r_{P,Q}:N_{\ca L}(P,Q)\to Hom_{\ca F}(P,Q)$ is the map which sends $g$ to the conjugation 
homomorphism $c_g:P\to Q$. We write $\r_P$ for the homomorphism 
$\r_{P,P}:Aut_{\ca T}(P)\to Aut_{\ca F}(P)$. 

Since $(S,\D,S)$ is a locality we have the category $\ca T_{\D}(S)$; and there is a functor 
$$ 
\e:\ca T_{\D}(S)\to\ca T 
$$ 
which is the identity map $\D\to\D$ on objects, and where $\e_{P,Q}$ is the inclusion map 
$N_S(P,Q)\to N_{\ca L}(P,Q)$. We write $\e_P$ for the inclusion $N_S(P)\to N_{\ca L}(P)$. Also, for 
$P,Q\in\D$ with $P\leq Q$, write $\iota_{P,Q}$ for $(\1)\e_{P,Q}$. 

\vskip .1in 
The definition from [BLO3] of a transporter system over a discrete $p$-toral group is 
embedded in the statement of the following result. We remind the reader that we always understand 
composition of morphisms in a category to be taken from left to right.

\proclaim {Proposition A.8} Let $(\ca L,\D,S)$ be a compact locality on $\ca F$ and let 
$$ 
\ca T_{\D}(S)@>\e>> \ca T @>\r>>\ca F 
$$  
be the pair of functors defined above. Then the following hold.  
\roster 

\item "{(A1)}" $\e$ is the identity on objects and $\r$ is the inclusion on objects. 

\item "{(A2)}" For each $P,Q\in\D$ the the group $Ker(\r_P)$ acts freely on $Mor_{\ca T}(P,Q)$ 
from the left (by composition), and $\r_{P,Q}$ is the orbit map for this action. Also, $Ker(\r_Q)$ 
acts freely on $Mor_{\ca T}(P,Q)$ from the right. 

\item "{(B)}" For each $P,Q\in\D$ the map $\e_{P,Q}$ is injective, and $\e_{P,Q}\circ\r_{P,Q}$ sends 
$g\in N_S(P,Q)$ to $c_g\in Hom_{\ca F}(P,Q)$. 

\item "{(C)}" For all $\phi\in Mor_{\ca T}(P,Q)$ and all $x\in P$, the following square commutes in $\ca T$. 
$$
\CD 
P                      @>\phi>>              Q   \\
@V{(x)\e_P}VV                       @VV{(x)(\r(\phi)\circ\e_Q)}V   \\ 
P                      @>\phi>>              Q 
\endCD   
$$

\item "{(I)}" Each $\ca F$-conjugacy class of subgroups in $\D$ contains a subgroup $P$ such that the 
image of $N_S(P)$ under $\e_P$ is a Sylow $p$-subgroup of $Aut_{\ca T}(P)$ (i.e. a subgroup of finite index 
relatively prime to $p$). 

\item "{(II)}" Let $\phi:P\to Q$ be a $\ca T$-isomorphism, and regard conjugation by $\phi$ as a 
mapping $c_\phi:Aut_{\ca T}(P)\to Aut_{\ca T}(Q)$. Let $P\norm\bar P\leq S$ and $Q\norm\bar Q\leq S$, 
be given, and suppose that $c_\phi$ maps $(\bar P)\e_P$ into $(\bar Q)\e_Q$. Then there exists 
$\bar\phi\in Mor_{\ca T}(\bar P,\bar Q)$ such that $\iota_{P,\bar P}\circ\bar\phi=\phi\circ\iota_{Q,\bar Q}$. 

\item "{(III)}" Let $P_1\leq P_2\leq P_3\leq\cdots$ be an increasing sequence of members of $\D$, and 
for each $i$ let $\psi_i:P_i\to S$ be a $\ca T$-homomorphism. Assume that 
$\psi_i=\iota_{P_i,P_{i+1}}\circ\psi_{i+1}$ for all $i$, and set $P=\bigcup\{P_i\}_{i=1}^\infty$. 
Then there exists $\psi\in Mor_{\ca T}(P,S)$ such that $\psi_i=\iota_{P_i,S}\circ\psi$ for all $i$. 

\endroster 
\endproclaim 

\demo {Proof} The condition (A1) is immediate from the definition of the functors $\e$ and $\r$. 
Under the identification of $Aut_{\ca T}(P)$ with $N_{\ca L}(P)$ we have composition in 
$Aut_{\ca T}(P)$ given by group multiplication in $N_{\ca L}(P)$, $Ker(\r_P)=C_{\ca L}(P)$, 
and similarly $Ker(\r_Q)=C_{\ca L}(Q)$. The actions defined in (A2) are then obviously free, and 
since $Aut_{\ca F}(P)\cong N_{\ca L}(P)/C_{\ca L}(P)$ we obtain the conclusion of (A2). 

The condition (B) is again immediate from the definition of $\e$. Now let $g\in N_{\ca L}(P,Q)$ and 
let $x\in P$. Regard $g$ as a $\ca T$-homomorphism $\phi:P\to Q$. Then $(x)\e_P\circ\phi$ is 
simply the product $xg$, while the composition $\phi\circ((x)(\r(\phi)\circ\e_Q))$ is the 
product $gx^g$. As $xg=gx^g$ we have the required commutativity of the diagram in (C). 

Each $\ca F$-conjugacy class of subgroups in $\D$ contains a subgroup $P$ such that $N_S(P)$ 
is a Sylow $p$-subgroup of $N_{\ca L}(P)$, by I.3.10. Thus (I) holds.  

Again let $g\in N_{\ca L}(P,Q)$ be a $\ca T$-isomorphism. Then $P^g=Q$, and $c_g$ is an isomorphism 
$N_{\ca L}(P)\to N_{\ca L}(Q)$. If $P\norm\bar P\leq S$ and $Q\norm\bar Q\leq S$ with 
$\bar P^g\leq\bar Q$, then $g\in N_{\ca L}(\bar P,\bar Q)$, and in this way $g$ is a 
$\ca T$-homomorphism $\bar\phi:\bar P\to\bar Q$. That is, (II) holds. 

Let $(P_i)_{i=1}^\infty$, $(\psi_i)_{i=1}^\infty$, and $P$ be given as in (III). Then $\psi_i$, written 
in full detail, is a triple $(g_i,P_i,S)$ where $P_i\leq S_{g_i}$. The ``inclusion morphism" 
$\iota_{P_i,P_{i+1}}$ is the triple $(\1,P_i,P_{i+1})$, and thus 
$$ 
\iota_{P_i,P_{i+1}}\circ\psi_{i+1}=(g_{i+1},P_i,S). 
$$ 
The hypothesis of (III) therefore translates into the statement that the sequence $(g_i)$ is a constant 
sequence $(g)$ where $P\leq S_g$. Taking $\psi=(g,P,S)$ then yields the conclusion of (III). 
\qed 
\enddemo

\definition {Definition A.9} Let $\ca F$ be a fusion system over the discrete $p$-toral group $S$ and 
let $\D$ be any $\ca F$-closed set of subgroups of $S$. A {\it transporter system} associated to $\ca F$ 
consists of a category $\ca T$ with $Ob(\ca T)=\D$, together with a pair of functors  
$$ 
\ca T_{\D}(S)@>\e>>\ca T @>\r>>\ca F, 
$$ 
which satisfy the conditions (A1), (A2), (B), (C), (I), (II), and (III) from the preceding proposition 
(where $\r_P:Aut_{\ca T}(P)\to Aut_{\ca F}(P)$ is an abbreviation for $\r_{P,P}$, and where 
$\iota_{P,Q}$ is an abbreviation for $(\1)\e_{P,Q}$ if $P,Q\in\D$ with $P\leq Q$). Write 
$\iota_P$ for $\iota_{P,P}$. 
\enddefinition

The definitions of linking system and of $p$-local compact group may now be given as follows, by 
[BLO3, Corollary A.5].

\definition {Definition A.10} The transporter system 
$$ 
\ca T_{\D}(S)@>\e>> \ca T @>\r>>\ca F 
$$  
is a {\it linking system} associated with $\ca F$ if the following conditions hold. 
\roster 

\item "{(1)}" $\ca F$ is order-saturated. 

\item "{(2)}" We have $P\in\D$ for each $\ca F$-centric subgroup $P\leq S$ such that 
$O_p(Out_{\ca F}(P))=1$. 

\item "{(3)}" For each $P$ in $\D$ the kernel of the homomorphism 
$\r_P:Aut_{\ca T}(P)\to Aut_{\ca F}(P)$ is discrete $p$-toral. 

\endroster 
In the special case that $\D$ is the set $\ca F^c$ of all $\ca F$-centric subgroups of $S$ we say 
that $(\e,\r)$ is a {\it $p$-local compact group}. 
\enddefinition

\proclaim {Proposition A.11} Let $(\ca L,\D,S)$ be a compact locality on $\ca F$, with 
$\D=\ca F^c$. Then the transporter system $(\e,\r)$ given by proposition A.8 is a $p$-local compact group. 
\endproclaim 

\demo {Proof} We need only verify the conditions (1) and (3) in the preceding definition, since the 
hypothesis that $\D=\ca F^c$ yields the remaining requirements. Condition (1) is given by A.7. 
Now let $P\in\D$. There is then an isomorphism $\a:N_{\ca L}(P)\to Aut_{\ca T}(P)$ given by 
$g\maps (g,P,P)$, and the kernel of $\r_P$ is then the image under $\a$ of $C_{\ca L}(P)$. 
As $\ca L$ is proper and $P\in\ca F^c$, 6.9 yields $C_{\ca L}(P)=Z(P)$, and so (3) holds. 
\qed 
\enddemo

Our goal now is to proceed in the opposite direction from that of the preceding result. Thus, starting with 
a $p$-local compact group, we aim now to construct a compact locality $(\ca L,\D,S)$.

\vskip .1in 
In what follows we fix the $p$-local finite group $(\ca T_{\D}(S)@>\e>>\ca T @>\r>>\ca F)$, with the 
abbreviations $\r_P$ and $\iota_{P,Q}$ as earlier. Write $\iota_P$ for the identity morphism $\iota_{P,P}$ 
in $Aut_{\ca T}(P)$. Condition A.9(B) implies that the image of $\iota_{P,P'}$ under $\r$ is the 
inclusion map $P\to P'$, so $\iota_{P,P'}$ is referred to as an {\it inclusion morphism} of $\ca T$. 
This leads to the following definition. 

\definition {Definition A.12} Let $P,Q,P',Q'\in\D$ with $P\leq P'$ and $Q\leq Q'$, and further let 
$\phi\in Mor_{\ca T}(P,Q)$ and $\phi'\in Mor_{\ca T}(P',Q')$. Then $\phi$ is an 
{\it extension} of $\phi'$, and $\phi'$ is a {\it restriction} of $\phi$ if  
$$ 
\iota_{P,P'}\circ\phi=\phi\circ\iota_{Q,Q'}. 
$$
\enddefinition

The following result collects the basic properties concerning the transporter system $(\e,\r)$.

\proclaim {Lemma A.13} 
\roster 

\item "{(a)}" Let $P,Q,R\in\D$, and let 
$$ 
P@>\bar{\phi}>>Q\quad\text{and}\quad Q@>\bar{\psi}>> R
$$ 
be $\ca F$-homomorphisms. Further, let $\psi\in Mor_{\ca T}(Q,R)$ with $\r(\psi)=\bar\psi$, and let 
$\l\in Mor_{\ca T}(P,R)$ with $\r(\l)=\bar\phi\circ\bar\psi$. Then there exists a unique 
$\phi\in Mor_{\ca T}(P,R)$ such that $\r(\phi)=\bar\phi$ and such that $\l=\phi\circ\psi$. 

\item "{(b)}" Let $\psi:P\to Q$ be a $\ca T$-morphism and let $P_0,Q_0\in\D$ with $P_0\leq P$ and 
with $Q_0\leq Q$. Suppose that $\r(\psi)$ maps $P_0$ into $Q_0$. There is then a unique 
$\ca T$-morphism $\psi_0:P_0\to Q_0$ such that $\psi$ is an extension of $\psi_0$. 

\item "{(c)}" A $\ca T$-homomorphism $\phi$ is a $\ca T$-isomorphism if and only if 
$\r(\phi)$ is an $\ca F$-isomorphism. 

\item "{(d)}" All morphisms of $\ca T$ are both monomorphisms and epimorphisms in the categorical 
sense. That is, we have left and right cancellation for morphisms in $\ca T$. 

\item "{(e)}" Let $\phi_0:P_0\to Q_0$ be a $\ca T$-morphism and let $P_0\leq P\leq S$ and 
$Q_0\leq Q\leq S$. Then there exists at most one extension of $\phi_0$ to a $\ca T$-homomorphism 
$P\to Q$. 

\item "{(f)}" Let $P,\bar P,Q,\bar Q$ be objects of $\ca T$, with $P\norm\bar P$ and with $Q\norm\bar Q$. 
Suppose that we are given a $\ca T$-isomorphism $\phi:P\to Q$ and an extension of $\phi$ to a 
$\ca T$-homomorphism $\bar\phi:\bar P\to \bar Q$. Then for each $x\in\bar P$ there is a 
commutative square: 
$$ 
\CD 
P                  @>\phi>>          Q \\ 
@V{x\d_{P,P}}VV                  @VV{y\d_{Q,Q}}V  \\ 
P                  @>\phi>>          Q 
\endCD 
$$ 
where $y$ is the image of $x$ under $\r(\bar\phi)$. 

\item "{(g)}" Every $\ca T$-morphism $\psi:P\to Q$ is the composite of a $\ca T$-isomorphism 
$\phi:P\to Q_0$ followed by an inclusion morphism $\iota_{Q_0,Q}$, where $Q_0$ is the image of $P$ 
under $\r(\psi)$.  

\endroster 
\endproclaim 

\demo {Proof} Points (a) through (d) constitute [BLO3, Proposition A.2], and (e) follows from the 
left cancellation in (d). For (f) one may appeal to the proof of [OV, Lemma 3.3(d)], as that proof 
does not depend on the finiteness of $S$. Finally, let $\phi:P\to Q$ be a $\ca T$-homomorphism and 
let $Q_0$ be the image of $P$ under $\r(\phi)$. Then there exists a restriction of $\phi$ to 
a $\ca T$-homomorphism $\phi_0:P\to Q_0$ by (b), and $\phi_0$ is then a $\ca T$-isomorphism by (c). 
This yields (g). 
\qed
\enddemo  

(A.14) By [BLO2, Section 3] it is a feature of an order-saturated fusion system $\ca F$ over a discrete 
$p$-toral group $S$ that there is a mapping $P\maps P^{\bull}$ from $Sub(S)$ into $Sub(S)$
having the following properties. 
\roster 

\item "{(1)}" $\{P^{\bull}\mid P\leq S\}$ is $\ca F$-invariant, and is 
the union of a finite number of $S$-conjugacy classes of subgroups of $S$. 

\item "{(2)}" For subgroups $P\leq Q\leq S$ we have $P^{\bull}\leq Q^{\bull}$ and 
$(P^{\bull})^{\bull}=P^{\bull}$. 

\item "{(3)}" For all $P,Q\leq S$ we have $N_S(P,Q)\sub N_S(P^{\bull},Q^{\bull})$. 

\item "{(4)}" For all $P,Q\leq S$, each $\ca F$-homomorphism $\a:P\to Q$ extends to an 
$\ca F$-homomorphism $\a^{\bull}:P^{\bull}\to Q^{\bull}$. 

\endroster 
In fact, we will not need (4) here. Rather, what we require is the following result concerning 
$p$-local finite groups.

\proclaim {Lemma A.15} Let $\ca T^{\bull}$ be the full subcategory of $\ca T$ whose set of objects is 
$\{P^{\bull}\mid P\in\D\}$. Then there is a functor 
$$ 
(-)^{\bull}:\ca T\to\ca T^{\bull}, 
$$ 
having the following properties. 
\roster 

\item "{(a)}" $(-)^{\bull}$ is the mapping $P\maps P^{\bull}$ on objects $P\in\D$. 

\item "{(b)}" $(-)^{\bull}$ restricts to the identity functor on $\ca T^{\bull}$. 

\item "{(c)}" For all $P,Q\in\D$ and all $\phi\in Mor_{\ca T}(P,Q)$, the image $\phi^{\bull}$ of 
$\phi$ under $(-)^{\bull}$ is an extension of $\phi$. 

\item "{(d)}" If $\a:X\to Y$ and $\phi:P\to Q$ are $\ca T$-morphisms such that $\phi$ is an 
extension of $\a$, then $\phi^{\bull}$ is an extension of $\a^{\bull}$. 

\endroster 
\endproclaim 

\demo {Proof} Points (a) through (c) are given by [JLL, Proposition 1.12]. By the same reference we 
have also the result that for all $X,P\in\D$ and all $g\in N_S(X,P)$ we have (in accord with A.14(3)) 
$(g)\e_{X,P})^{\bull}=(g)\e_{X^{\bull},P^{\bull}}$. In particular, by taking $X\leq P$ and $g=\1$ we obtain 
$(\iota_{X,P})^{\bull}=\iota_{X^{\bull},P^{\bull}}$. If $\a$ and $\phi$ are given as in (d), so that 
$\iota_{X,P}\circ\phi=\a\circ\iota_{Y,Q}$, the functoriality of $(-)^{\bull}$ now yields 
$\iota_{X^{\bull},P^{\bull}}\circ\phi^{\bull}=\a^{\bull}\circ\iota_{Y^{\bull},Q^{\bull}}$. 
Thus (d) holds. 
\qed 
\enddemo

\proclaim {Lemma A.16} Let $\phi_0:P_0\to Q_0$, $\phi:P\to Q$, and $\phi':P'\to Q'$ be 
$\ca T$-isomorphisms, and suppose that both $\phi$ and $\phi'$ are extensions of $\phi_0$. 
\roster 

\item "{(a)}" If $P=P'$ or if $Q=Q'$, then $\phi=\phi'$. 

\item "{(b)}" There is a unique extension of $\phi_0$ to an isomorphism $\psi:P\cap P'\to Q\cap Q'$, 
and both $\phi$ and $\phi'$ are extensions of $\psi$. 

\endroster 
\endproclaim 

\demo {Proof} Assume that (a) is false. We may take $P=P'$, since the case where $Q=Q'$ will then follow 
by considering the inverses of the given $\ca T$-isomorphisms. Note that if $\phi^{\bull}=(\phi')^{\bull}$ 
then $\phi=\phi'$ by restriction. Since $\phi^{\bull}$ and $(\phi')^{\bull}$ are extensions of 
$(\phi_0)^{\bull}$ by A.15(b), it therefore suffices to consider the case where $\phi_0$, $\phi$, and 
$\psi$ are $\ca T^{\bull}$-isomorphisms. The finiteness condition A.14(1) then yields the existence 
of a counter-example $(\phi_0,\phi,\phi')$ to (a) in which $|P_0|$ is maximal. 

Let $x\in N_P(P_0)$, let $y$ be the image of $x$ under $\r(\phi)$, and let $y'$ be the image of $x$ 
under $\r(\phi')$. We appeal to A.14(f) with $(P_0,Q_0,N_P(P_0),N_Q(Q_0))$ in the role of  
$(P,Q,\bar P,\bar Q)$, and obtain 
$$ 
\phi_0\i\circ(x)\e_{P_0,P_0}\circ\phi_0=(y)\e_{Q_0,Q_0}=(y')\e_{Q_0,Q_0}.  
$$ 
As $\e_{Q_0,Q_0}$ is injective (by A.9(B)) we get $y=y'$, and thus $\r(\phi)$ and $\r(\phi')$ 
agree on $P_1:=N_P(P_0)$. Let $Q_1$ be the image of $P_1$ under $\r(\phi)$. By A.13(b) there is a 
restriction $\phi_1:P_1\to Q_1$ of $\phi$ and a restriction $\phi_1':P_1\to Q_1$ of 
$\phi'$, and then $\phi_1=\phi_1'$ by A.13(e). Now $(\phi_1,\phi,\phi')$ is a counter-example to (a)   
with $|P_1|>|P_0|$, in violation of the maximality of $|P_0|$. This contradiction 
completes the proof of (a). 

Set $X=P\cap P'$ and $Y=Q\cap Q'$. Then $\phi$ and $\phi'$ have 
restrictions $\psi:X\to(X)(\r(\phi))$ and $\psi':X\to(X)(\r(\phi'))$ which, in turn, restrict to 
$\phi_0$. Then (a) yields $\psi=\psi'$, and this establishes (b). 
\qed 
\enddemo 

Define a relation $\up$ on the set $Iso(\ca T)$ $\ca T$-isomorphisms by $\phi\up\phi'$ if 
$\phi'$ is an extension of $\phi$. We may also write $\phi'\down\phi$ to indicate that $\phi$ is a 
restriction of $\phi'$. 

\proclaim {Lemma A.17} The following hold.
\roster

\item "{(a)}" The relation $\up$ induces a partial order on $Iso(\ca T)$.  

\item "{(b)}" The relation $\up$ respects composition of morphisms. That is, if $\phi\up\phi'$ and 
$\psi\up\psi'$, and the compositions $\phi\circ\psi$ and $\phi'\circ\psi'$ are defined, then 
$(\phi\circ\psi)\up(\phi'\circ\psi')$. 

\item "{(c)}" For each $\ca T$-isomorphism $\a$ there exists a unique $\ca T$-isomorphism $\phi$ 
such that $\phi$ is maximal with respect to $\up$ and such that $\a\up\phi$. 

\endroster 
\endproclaim 

\demo {Proof} For points (a) and (b) we repeat the proof of [Che1, Lemma X.7]. The transitivity of  
$\up$ is immediate. Suppose that both $\phi\up\phi'$ and $\phi\down\phi'$, where $\phi\in Iso_{\ca T}(P,Q)$ 
and $\phi'\in Iso_{\ca T}(P',Q')$. Then $P=P'$, $Q=Q'$, $\iota_{P,P'}=\iota_P$, and $\iota_{Q,Q'}=\iota_Q$. 
Further, $\iota_P\phi'=\phi\circ\iota_{Q}$ and then $\phi'=\phi$ since $\iota_P$ and $\iota_Q$ are identity 
morphisms in $\ca T$. Thus (a) holds. 

Suppose that we are given $\phi\up\phi'$ and $\psi\up\psi'$, with $\phi\circ\psi$ and $\phi'\circ\psi'$ 
defined on objects $P$ and $P'$ respectively. Set $Q=P\phi$ and $R=Q\psi$, and set $Q'=P'\phi'$ and 
$R'=Q'\psi'$. The following diagram, in which the vertical arrows are inclusion morphisms, demonstrates that 
$\phi\circ\psi\up\phi'\circ\psi'$. 
$$
\CD
P'         @>\phi'>>    Q'       @>\psi'>>       R' \\ 
@AAA                   @AAA                     @AAA  \\ 
P          @>>\phi>     Q        @>>\psi>        R 
\endCD 
$$ 
This yields (b). 

Let $\a\in Iso(\ca T)$. The finiteness condition A.14(1), together with A.15(c,d), yields the existence 
of at least one $\ca T$-isomorphism $\phi$ such that $\a\up\phi$ and such that $\phi$ is maximal with 
respect to $\up$. Assuming now that $\a$ is a counter-example to (c), there then exists an 
$\up$-maximal $\ca T$-isomorphism $\phi'$ with $\a\up\phi'$ and with $\phi\neq\phi'$. Write 
$\a:X\to Y$, $\phi:P\to Q$, and $\phi':P'\to Q'$. We may again apply A.14(1) in conjunction with 
A.15(c,d) in order to obtain such a triple $(\a,\phi,\phi')$ in which $|X|$ has been maximized. 

Set $P_1=N_P(X)$ and $Q_1=N_{Q}(Y)$, and similarly define $P_1'$ and  $Q_1'$. Set 
$X_1=\<P_1,P_1'\>$ and $Y_1=\<Q_1,Q_1'\>$. Let $\l:Aut_{\ca T}(P_0)\to Aut_{\ca T}(Q_0)$ be the 
isomorphism induced by conjugation by $\a$. Then A.13(f) implies that $\l$ maps 
$(X_1)\e_{P_0}$ to $(Y_1)\e_{Q_0}$. Condition (II) in the definition of transporter system 
then yields the existence of an extension of $\a$ to an isomorphism $\a_1:X_1\to Y_1$. 
Let $\phi_1$ be the restriction of $\phi$ to an isomorphism $P_1\to Q_1$. Then 
$\phi\down\phi_1\up\a_1$. Let $\psi$ be an $\up$-maximal extension of $\a_1$. If $\phi=\psi$ 
then $X_1=X$, whence $P=X=P'$, and then A.16(a) yields $\phi=\phi'$. Thus $\phi\neq\psi$, so 
$(\phi_1,\phi,\psi)$ provides a counter-example to (c). The maximality of $|X|$ then yields 
$X=N_P(X)$, and we again obtain $X=P$. Then $\phi=\a$, $\a$ is $\up$-maximal, and again $\a=\phi'$. 
Thus we have a contradiction, proving (c). 
\qed 
\enddemo

Let $\equiv$ be the equivalence relation on $Iso(\ca T)$ generated by $\up$, and let 
$$ 
\ca L=Iso(\ca T)/\equiv 
$$ 
be the set of equivalence classes. For $\phi\in Iso(\ca T)$ we write $[\phi]$ for the equivalence 
class containing $\phi$.

\proclaim {Lemma A.18} Let $f\in\ca L$. 
\roster 

\item "{(a)}" There is a unique $\phi\in f$ such that $\psi$ is $\up$-maximal in the poset $Iso(\ca T)$. 
Moreover, we then have $\a\up\phi$ for all $\a\in f$, and $\phi\i$ is the unique $\up$-maximal member 
of $[\phi\i]$. 

\item "{(b)}" The unique maximal $\phi\in f$ is a $\ca T^{\bull}$-isomorphism. 

\item "{(c)}" $f\cap Iso_{\ca T}(P,Q)$ has cardinality at most 1 for any 
$P,Q\in Ob(\ca T)$. 

\endroster 
\endproclaim 

\demo {Proof} There exists at least one $\up$-maximal member $\phi\in f$ by A.14(1) with A.15(c,d). Suppose 
that $\phi$ and $\phi'$ are distinct $\up$-maximal members of $f$. 
As $\phi\equiv\phi'$ there is a sequence $\s=(\psi_0,\cdots,\psi_n)$ of members of $f$ with $\phi=\psi_0$, 
$\phi'=\psi_n$, and such that for all $i$ with $1\leq i\leq n$ we have either $\psi_{i-1}\up\psi_i$ or 
$\psi_{i-1}\down\psi_i$. Assume that the pair $(\phi,\phi')$ has been chosen so that $n$ is as small as 
possible. The maximality of $\phi$ implies $\phi\down\psi_1$ and that $\psi_2$ is not a restriction of 
$\phi$. Since $\psi_2$ is the restriction of some maximal isomorphism (necessarily in $f$), we obtain 
$n=2$. Thus $\phi\down\psi_1\up\phi'$. As this violates A.17(c) we obtain the uniqueness asserted in (a). 
Moreover, A.17(c) then implies that each $\a\in f$ extends to $\phi$. The inverse of any extension of 
$\a$ is an extension of $\a\i$, so (a) holds. Point (b) then follows from A.15(c).  

In order to prove (c), let $\psi,\psi'\in f\cap Iso_{\ca T}(P,Q)$. Then both $\psi$ and $\psi'$ are 
restrictions of a single $\phi\in f$, by (a). Now A.13(b) yields $\psi=\psi'$. 
\qed 
\enddemo

Define $\bold D$ to be the set of words $w=(f_1,\cdots,f_n)\in\bold W(\ca L)$ such that there exists a 
sequence $(\phi_1,\cdots,\phi_n)$ of $\ca T$-isomorphisms with $\phi_i\in f_i$, and a sequence 
$(P_0,\cdots,P_n)$ of members of $\D$, such that each $\phi_i$ is a $\ca T$-isomorphism $P_{i-1}\to P_i$.  
As in section I.2 we may say also that $w\in\bold D$ via $(P_0,\cdots,P_n)$, or via $P_0$. Define 
$$
\Pi:\bold D\to\ca L 
$$
by $\Pi(w)=f$, where $f$ is the unique maximal element of $[\phi_1\circ\cdots\circ\phi_n]$ given by 
A.18(a). That $\Pi$ is well-defined follows from A.17(b). Set 
$\1=[\iota_S]$, and for any $f\in\ca L$ let $f\i$ be the equivalence class of 
$\phi\i$, where $\phi$ is the unique maximal member of $f$.

\proclaim {Proposition A.19} $\ca L$ with the above structures is a partial group. Moreover, the 
following hold. 
\roster 

\item "{(a)}" For any $x\in S$, the $\equiv$-class $[(x)\e_S]$ is the set of all $(x)\e_{P,Q}$ such that 
$P^x=Q$, and $(x)\e_S$ is the maximal member of $[(x)\e_S]$. 

\item "{(b)}" $[\iota_S]=\{\iota_P\mid P\in Ob(\ca T)\}$, and $\iota_S$ is the maximal member of its class. 

\item "{(c)}" For any $\phi\in Iso(\ca T)$, $[\phi\i]$ is the set of 
inverses of the members of $[\phi]$. 

\endroster 
\endproclaim 

\demo {Proof} We first check via definition I.1.1 that $\ca L$ is a partial group. For any $f\in\ca L$ 
the word $(f)$ of length 1 is in $\bold D$ since $f$ is represented by a $\ca T$-isomorphism. If $w\in\bold D$ 
and $w=u\circ v$ then it is immediate from the definition of $\bold D$ that both $u$ and $v$ are in 
$\bold D$. Thus the condition I.1.1(1) in the definition of partial group is satisfied. By definition of 
$\Pi$ we have $\Pi(f)=f$ for $f\in\ca L$, so I.1.1(2) holds. Condition I.1.1(3) is a straightforward 
consequence of associativity of composition of isomorphisms in $\ca T$. 

That the inversion map $f\maps f\i$ is an involutory bijection follows from 
A.18(a). Now let $u=(f_1,\cdots,f_n)\in\bold D$ via $(P_0,\cdots,P_n)$, and 
set  $u\i=(f_n\i,\cdots,f_1\i)$. Then $u\i\in\bold D$ via $(P_n,\cdots,P_0)$, 
so $u\i\circ u\in\bold D$. One obtains a representative in the class 
$\Pi(u\i\circ u)$ via a sequence of cancellations $\phi_k\i\circ\phi=\iota_{P_k}$, for 
representatives $\phi_k\in f_i$, so $\Pi(u\i\circ u)$ is the equivalence 
class containing $\iota_{P_0}$. Since $\iota_{P_0}\up\iota_S$, and since 
$\1=[\iota_S]$ by definition, we get $\Pi(u\i\circ u)=\1$. Thus I.1.1(4) holds 
in $\ca L$, and $\ca L$ is a partial group. 

We now prove (a). Let $P\leq P'$ and $Q\leq Q'$ in $\D$, and let $x$ be an element of $S$ such that 
$P^x=Q$ and $(P')^x=Q'$. The functoriality of $\e$ yields 
$$  
(\1)\e_{P,P'}\circ(x)\e_{P',Q'}=(x)\e_{P,Q'}=(x)\e_{P,Q}\circ(\1)\e_{Q,Q'}, 
$$ 
which means that $(x)\e_{P,Q}(g)\up(x)\e_{P',Q'}$. In particular, we get 
$(x)\e_{P,Q}(g)\up(x)\e_S$. In order to complete the proof of (a), it now 
suffices to show that for any $\phi\in Iso_{\ca T}(P,Q)$ with $(x)\e_S\equiv\phi$, we have  
$\phi=(x)\e_{P,Q}$. 
Suppose false, and let $\s=(\phi_1,\cdots,\phi_n)$ be a sequence of 
$\ca T$-isomorphisms with $\phi=\phi_1$, $(x)\e_S=\phi_n$, and with either 
$\phi_i\up\phi_{i+1}$ or $\phi_i\downarrow\phi_{i+1}$ for all $i$ with 
$1\leq i<n$. Among all $(\phi,P,Q)$ with $\phi\neq(x)\e_{P,Q}$ and 
$(x)\e_S\equiv\phi$, choose $(\phi,P,Q)$ so that the length of such a chain 
$\s$ is as small as possible. Set $\psi=\phi_2$. Then $\psi=(x)\e_{X,Y}$, 
where $X$ and $Y$ are objects of $\ca T$ with $X^x=Y$. Suppose 
$\phi\uparrow\psi$. Applying the functor $\r$ to the commutative diagram 
$$ 
\CD 
     X     @>(x)\e_{X,Y}>>    Y  \\
@V\iota_{P,X}VV               @VV\iota_{Q,Y}V \\ 
     P    @>\phi>>            Q
\endCD,
$$
and applying condition (B) in the definition of transporter system to $\r(\e_{X,Y}(g))$, we conclude that 
$\r(\phi)$ is the restriction of $c_g$ to the homomorphism $\r(\phi):P\to Q$. In particular, we get $P^x=Q$, 
so that also $(x)\e_{P,Q}$ is a restriction of $(x)\e_{X,Y}$. Then $\phi=(x)\e_{P,Q}$ by A.13(d), and 
contrary to hypothesis. On the other hand, if $\phi\down\psi$, then $\phi=(x)\e_{P,Q}$ by A.16(a), again 
contrary to hypothesis. This completes the proof of (a), and then (b) is the special case of (a) given
by $x=1$. 

Let $f=[\phi]$ be an equivalence class, with $\phi$ maximal in $f$. One 
checks (by reversing pairs of arrows in the appropriate diagrams) that if 
$\psi$ is a $\ca T$-isomorphism, and $\psi$ is a restriction of $\phi$, then 
the $\ca T$-isomorphism $\psi\i$ is a restriction of $\phi\i$. Point (c) 
follows from this observation. 
\qed 
\enddemo

\definition {Remark} In view of A.19(a) there will be no harm in writing $x$ to denote the equivalence 
class $[(x)\e_S]$, for $x\in S$. That is to say that from now on we shall identify $S$ with the image 
of $S$ under the composition of $\e_S$ with the projection $Iso(\ca T)\to\ca L$. 
\enddefinition

\proclaim {Lemma A.20} Let $\phi:Z\to W$ be a $\ca T$-isomorphism, maximal in 
its $\equiv$-class. Let $X$ and $Y$ be objects of $\ca T$ contained in $Z$, 
and let $U$ and $V$ be the images of $X$ and $Y$, respectively, under 
$\r(\phi)$. Suppose that there exist elements $x$ and $x'$ in $S$ such that 
the following diagram commutes. 
$$
\CD 
X                  @>{\phi\mid_{X,U}}>>            U   \\ 
@V{(x)\e_{X,Y}}VV                         @VV{(x')\e_{U,V}}V   \\ 
Y                  @>{\phi\mid_{Y,V}}>>            V 
\endCD\tag* 
$$ 
Then $x\in Z$, and $x'$ is the image of $x$ under $\r(\phi)$. 
\endproclaim

\demo {Proof} Let $\phi'$ be the composition (in right-hand notation) 
$$
\phi'=(x\i)\e_{Z^x,Z}\circ\phi\circ(x')\e_{W,W^{x'}}.  
$$
Thus, $\phi'\in Iso_{\ca T}(Z^x,W^{x'})$, and the commutativity of (*) 
yields $\phi\equiv\phi'$. The maximality 
of $\phi:Z\to W$ implies that $Z^x\leq Z$ and $W^{x'}\leq W$. That is, 
$x\in N_S(Z)$ and $x'\in N_S(W)$. There is then a commutative diagram as follows. 
$$
\CD
Z                 @>{\phi}>>            W   \\
@V{(x)\e_{Z}}VV                   @VV{(x')\e_{W}}V   \\
Z                 @>{\phi}>>            W
\endCD
$$
Condition (II) in the definition of transporter system implies that there is an extension of $\phi$ to a 
$\ca T$-isomorphism $\<Z,x\>\to\<W,x'\>$, and the maximality of $\phi$ then yields $x\in Z$ and $x'\in W$. 
Condition (C) in the definition of transporter system implies that $x'$ is the image under $\r(\phi)$ of $g$. 
\qed
\enddemo

\proclaim {Corollary A.21} Let $f\in\ca L$ and let $P\in\D$ with the property 
that, for all $x\in P$, $(f\i,x,f)\in\bold D$ and $\Pi(f\i,x,f)\in S$. Let 
$Q$ be the set of all such products $\Pi(f\i,x,f)$. Then $Q\in\D$ and there 
exists $\psi\in f$ such that $\psi\in Iso_{\ca T}(P,Q)$. 
\endproclaim 

\demo {Proof} As $(f\i,x,f)\in\bold D$ there exist $U,X,Y,V\in\D$ and 
representatives $\psi$ and $\bar\psi$ of $f$ such that 
$$ 
\CD 
U@>\bar\psi\i>>X@>(x)\e_{X,Y}>>Y@>\psi>>V 
\endCD 
$$ 
is a chain of $\ca T$-isomorphisms, and where the middle arrow in the 
diagram is given by A.19(a). As $\Pi(f\i,x,f)\in S$ there exists $x'\in S$ such that 
$\bar\psi\i\circ(x)\e_{X,Y}\circ\psi=(x')\e_{U,V}$. Let $\phi:Z\to W$ be the 
maximal element of $f$. Then A.20 implies that $x\in Z$, and 
$x'$ is the image of $x$ under $\r(\phi)$. In particular, we 
have $P\leq Z$ and $Q\leq W$, and we may therefore take $X=Y=P$ 
and $U=V=Q$, obtaining $\psi\in Iso_{\ca T}(P,Q)$. 
\qed  
\enddemo 

\proclaim {Lemma A.22} Let $\psi:P\to Q$ be a $\ca T$-isomorphism, and let $f=[\psi]$ be the 
equivalence class of $\psi$. Then $(f\i,x,f)\in\bold D$ for all $x\in P$, and $P^f=Q$ in the partial 
group $\ca L$. Moreover, the conjugation map $(f\i,x,f)\maps\Pi(f\i,x,f)$ is equal to $\r(\psi)$. 
\endproclaim

\demo {Proof} For any $x\in P$, we have the composable sequence 
$$
Q@>\psi\i>>P@>(x)\e_P>>P@>\psi>>Q  
$$
of $\ca T$-isomorphisms, so $(f\i,x,f)$ is in $\bold D$. By Condition (C) in the definition of transporter 
system then yields $\psi\i\circ(x)\e_P\circ\psi=(x')\e_Q$, where $x'=(x)(\r(\psi))\in Q$. 
The class $[(x')\e_Q]$ is the same as $[(x')\e_S]$ by 
A.19(a); and we recall that we have introduced the convention to denote this 
class simply as $x'$. Thus $x^f=x'$, and so $P^f\sub Q$. Similarly $Q^{f\i}\sub P$, from which one 
deduces that the conjugation map $c_f:P\to Q$ is surjective. Injectivity of $c_f$ follows from 
left and right cancellation in the partial group $\ca L$, so $P^f=Q$. The final assertion of the 
lemma is given by the observation, made above, that $x'=(x)(\r(\psi))$. 
\qed 
\enddemo

\proclaim {Theorem A.23} Let $(\ca T_{\D}(S)@>\e>> \ca T @>\r>>\ca F)$ be a $p$-local compact group, 
and let $\ca L=Iso(\ca T)/\equiv$ be the partial group given by A.19. For $x\in S$, identify $x$ with 
the $\equiv$-class of the $\ca T$-automorphism $(x)\e_S$ of $S$. Then $(\ca L,\D,S)$ is a 
compact locality on $\ca F$. 
\endproclaim 

\demo {Proof} We have seen in A.19 that $\ca L$ is a partial group, and the remark following A.19 
shows how to identify $S$ with its image under the composition of $\e_S$ with the quotient map 
$Iso(\ca T)\to\ca L$. In order to show that $(\ca L,\D)$ 
is objective, begin with $w=(f_1,\cdots,f_n)\in\bold D$. By definition, there exist representatives 
$\psi_i$ of the classes $f_i$, and a sequence $(P_0,\cdots,P_n)$ of objects in $\D$, such that each 
$\psi_i$ is a $\ca T$-isomorphism $P_{i-1}\to P_i$. Then $P_{i-1}^{f_i}=P_i$ for all $i$, 
by A.22. Conversely, given $w=(f_1,\cdots,f_n)\in\bold W(\ca L)$, and given 
$(P_0,\cdots,P_n)\in\bold W(\D)$ with $P_{i-1}^{f_i}=P_i$ for all $i$, it follows from A.21 that 
$w\in\bold D$. Thus, $(\ca L,\D)$ satisfies the condition (O1) in the definition I.2.1 of 
objective partial group. Since $\D=\ca F^c$ is $\ca F$-closed we also have (O2), and thus $(\ca L,\D)$ 
is objective. As $\D$ is a set of subgroups of $S$, $(\ca L,\D,S)$ is then a pre-locality (as defined 
in I.2.6). 

The conjugation maps $c_g:S_g\to S$ for $g\in\ca L$ are $\ca F$-homomorphisms, by A.21 and A.9(C). Thus 
$\ca F_S(\ca L)$ is a subsystem of $\ca F$. Assuming now that $\ca F\neq\ca F_S(\ca L)$, there exists 
an $\ca F$-isomorphism $\b:X\to Y$ such that $\b$ is not an $\ca F_S(\ca L)$-homomorphism. By A.14(4) 
we may take $X=X^{\bull}$ and $Y=Y^{\bull}$, and by the finiteness condition in A.14(1) we may then 
assume that from among all $\ca F$-isomorphisms which are not $\ca F_S(\ca L)$-homomorphisms, 
$\b$ has been chosen so that $|X|$ is as large as possible. If $X\in\D$ then $Y\in\D$ and the surjectivity 
of $\r_{X,Y}$ (condition (A2) in definition A.9) implies that $\b$ is an $\ca F_S(\ca L)$-homomorphism, 
so in fact $X\notin\D$. In particular $X<S$, and so $X<N_S(X)$. Similarly $Y<N_S(Y)$.  

As $\ca F$ is order-saturated there exists a fully order-normalized $\ca F$-conjugate $Z$ of $X$, 
and there then exist $\ca F$-homomorphisms $\eta_1:N_S(X)\to N_S(Z)$ 
and $\eta_2:N_S(Y)\to N_S(Z)$ such that $X\eta_1=Z=Y\eta_2$. Each $\eta_i$ is an 
$\ca F_S(\ca L)$-homomorphism by the maximality in the choice of $X$, and it then suffices to show 
that the $\ca F$-automorphism $\a=\eta_1\i\circ\b\circ\eta_2$ of $Z$ is an $\ca F_S(\ca L)$-homomorphism. 
As $\ca F$ is order-receptive, $\a$ extends to an $\ca F$-automorphism $\bar{\a}$ of 
$C_S(Z)Z$. But $C_S(Z)Z$ is centric in $\ca F$, so $C_S(Z)Z\in\D$, and $\bar{\a}$ is then 
an $\ca F_S(\ca L)$-homomorphism. The same is then true of $\a$, and so we have shown that 
$\ca F=\ca F_S(\ca L)$. 

Set $\G=\{P^{\bull}\mid P\leq S\}$. For $P,Q\in\G$ we have (by A.14(2)) 
$$ 
(P\cap Q)^{\bull}\leq P\cap Q^{\bull}=P\cap Q, 
$$ 
and thus $\G$ is closed under finite intersections. Let $g\in\ca L$ and let $\psi:P\to Q$ be the unique 
$\up$-maximal representative of $g$. Then $P=S_g$ by A.22, and then $S_g\in\G$ by A.15(c). We may 
prove by induction on the length of $w\in\bold W(\ca L)$ that $S_w\in\G$. Namely, write 
$w=(g)\circ v$ where $g\in\ca L$ and $v\in\bold W(\ca L)$ with $S_v\in\G$. Then 
$S_w=(S_{g\i}\cap S_v)^{g\i}\in\G$ since, as we have seen, $\G$ is closed with respect to finite 
intersections and since (by A.14(1)) $\G$ is $\ca F$-invariant. The finiteness condition in 
A.14(1) now implies that $\bigcap\{S_w\mid w\in \bold X\}\in\G$ for each non-empty subset $\bold X$ 
of $\bold W(\ca L)$. This shows that the poset $\Omega_S(\ca L)$ defined in I.2.11 is finite-dimensional. 

Let $P\in\D$ and let $\a_P:Aut_{\ca T}(P)\to N_{\ca L}(P)$ be the mapping $\psi\maps[\psi]$. Then 
$\a_P$ is a homomorphism by A.17(b), $\a_P$ is injective by A.19(c), and $\a_P$ is surjective by A.20. 
 Thus $\a_P$ is an isomorphism, and then condition (I) in definition A.9 implies that $N_{\ca L}(P)$ is 
virtually $p$-toral. As $\Omega_S(\ca L)$ is finite-dimensional, all subgroups of $\ca L$ are then virtually 
$p$-toral by I.2.17. In particular, if $\bar S$ is a $p$-subgroup of $\ca L$ containing $S$ then $\bar S$ 
is discrete $p$-toral. If $S<\bar S$ then A.2(c) yields $S<N_{\bar S}(S)$, which is contrary to condition (I) 
in A.5. Thus $S$ is a maximal $p$-subgroup of $\ca L$, and we have established that $\ca L$ is a  
locality on $\ca F$. Notice that A.22 implies that $\a_P$ restricts to an isomorphism 
$Ker(\r_P)\to C_{\ca L}(P)$. As $(\e,\r)$ is a $p$-local compact group, $Ker(\r_P)$ is a $p$-group, 
and thus $N_{\ca L}(P)$ is of characteristic $p$. That is, $\ca L$ satisfies the condition (PL2) in 
the definition (6.7) of proper locality. Condition (PL1), that $\ca F^{cr}$ be contained in $\D$, is 
given by $\D=\ca F^c$. Condition (PL3), that $S$ has the normalizer-increasing property, is given by 
A.2(c). Thus $\ca L$ is proper, and the proof is complete. 
\qed 
\enddemo

\proclaim {Theorem A.24} Let $(\ca L,\D,S)$ be a compact locality on $\ca F$, such that $\D$ 
is the set $\ca F^c$ of $\ca F$-centric subgroups of $S$. Let $(\ca T_{\D}(S)@>\e>> \ca T @>\r>>\ca F)$ 
be the $p$-local compact group constructed from $\ca L$ as in A.8 and A.11, and let $(\ca L',\D,S)$ be 
the compact locality constructed from $(\e,\r)$ as in A.19 and A.23. Then the mapping 
$$ 
\Phi:\ca L\to\ca L', 
$$  
which sends $g\in\ca L$ to the $\equiv$-class of the $\ca T$-isomorphism $(g,S_g,S_{g\i})$ (and with 
the identifications given by the remark following A.19) is an isomorphism of partial groups which restricts 
to the identity map on $S$. 
\endproclaim 

\demo {Proof} Let $\Phi^*:\bold W(\ca L)\to\bold W(\ca L')$ be the mapping induced by $\Phi$, and let 
$w=(g_1,\cdots,g_n)\in\bold D(\ca L)$ via $(P_0,\cdots,P_n)$, with $P_0=S_w$. We shall denote the 
$\equiv$-class of a $\ca T$-isomorphism $(g,P,Q)$ by $[g,P,Q]$. Then 
$w\Phi^*=([g_1,P_0,P_1],\cdots,[g_n,P_{n-1},P_n])$, and $w\Phi^*\in\bold D(\ca L')$ via $(P_0,\cdots,P_n)$ 
by A.22. The definition of the product $\Pi'$ in $\ca L'$ then yields 
$$ 
\Pi'(w\Phi^*)=[\pi(w),P_0,P_n]=(\Pi(w))\Phi, 
$$ 
and thus $\Phi$ is a homomorphism of partial groups. 

Recall that for $P,P'\in\D$ with $P\leq P'$, we have $\iota_{P,P'}=(\1,P,P')$. It follows that the extensions 
of a $\ca T$-isomorphism $(f,P,Q)$ are of the form $(f,P',Q')$, and this implies that $\Phi$ is injective. 
Since $\ca L'=Im(\Phi)$ by definition, $\Phi$ is a bijection, and we now leave it to the reader to verify 
that $\Phi\i$ is a homomorphism. For $x\in S$ we have identified $x$ with $[x,S,S]$, so $\Phi$ restricts to 
the identity map on $S$. 
\qed
\enddemo

\proclaim {Theorem A.25} Let $(\ca L,\D,S)$ and $(\ca L',\D,S)$ be compact localities on $\ca F$, having 
the same set of objects. Then there exists an isomorphism $\a:\ca L\to\ca L'$ of partial groups, 
such that $P\a=P$ for all $P\in\D\cap\ca F^c$. In particular, $\a$ restricts to an automorphism of $S$. 
\endproclaim 

\demo {Proof} By Theorem A1 we may assume without loss of generality that $\D$ is equal to the set $\ca F^c$ 
of $\ca F$-centric subgroups of $S$. Let $(\ca T_{\D}(S)@>\e>>\ca T@>\r>>\ca F)$ be the $p$-local 
compact group constructed from $(\ca L,\D,S)$ via A.8 and A.11, and let 
$(\ca T_{\D}(S)@>\e'>>\ca T'@>\r'>>\ca F)$ be the $p$-local compact group similarly constructed from 
$(\ca L',\D,S)$. 

Following [BLO2] (but with notation which reflects our preference for right-hand composition of morphisms), 
we define the {\it orbit category} $\ca O=\ca O^c(\ca F)$ to be the category whose set of 
objects is $\D$, with 
$$ 
Mor_{\ca O}(P,Q)=Hom_{\ca F}(P,Q)/Inn(Q). 
$$ 
That is, the $\ca O$-morphisms $P\to Q$ are the sets 
$$ 
[\phi]=\{\phi\circ c_x\mid x\in Q\}, 
$$ 
where $\phi$ is an $\ca F$-homomorphism $P\to Q$. If also $\psi:Q\to R$ is an $\ca F$-homomorphism then 
one has the well-defined composition $[\phi]\circ[\psi]=[\phi\circ\psi]$. 

There is a contravariant functor 
$$ 
\ca Z:\ca O^{op}\to Ab
$$ 
(where $Ab$ is the category of abelian groups), given by $\ca Z(P)=Z(P)$ on objects, and defined in the 
following way on $\ca O$-morphisms. If  $\phi:P\to Q$ is an $\ca F$-homomorphism (with $P$ and $Q$ 
centric in $\ca F$), then $\ca Z$ sends the $\ca O$-homomorphism $[\phi]$ to the homomorphism 
$Z(Q)\to Z(P)$ obtained as the composition of the inclusion map $Z(Q)\to Z(P\phi)$ followed by 
the map $\phi_0\i:Z(P\phi)\to Z(P)$, where $\phi_0$ is the $\ca F$-isomorphism $P\to P\phi$ induced 
by $\phi$. The main result of [LL] is: If $\ca F$ is saturated then the higher limit functors 
${\underset\leftarrow\to{lim}}^k(\ca Z)$ are trivial for all $k\geq 2$. 

There is a direct analogy with the theory of group extensions, which estabishes that  
the vanishing of ${\underset\leftarrow\to{lim}}^3(\ca Z)$ implies the existence of a 
$p$-local compact group $(\ca T_{\D}(S)@>\e>> \ca T @>\r>>\ca F)$. Such a $p$-local compact group 
can be viewed as an ``extension" of $\ca O$, in the sense that there is a functor 
$$ 
\ca T@>\s>>\ca O, 
$$ 
such that $\s$ induces the identity map $Ob(\ca T)\to Ob(\ca O)$ (i.e. the identity map on $\D$), and 
such that the image of a $\ca T$-morphism $\psi:P\to Q$ under $\s$ is equal to $[\phi]$, where 
$\phi$ is the $\ca F$-homomorphism $(\psi)\r$. The vanishing of ${\underset\leftarrow\to{lim}}^2(\ca Z)$ 
yields the uniqueness of this extension, up to isomorphism. That is, if 
$(\ca T_{\D}(S)@>\e'>> \ca T' @>\r'>>\ca F)$ is another $p$-local compact group, then there is 
an isomorphism $\b:\ca T\to\ca T'$ of categories, such that the following diagram commutes:  
$$ 
\CD
\ca T    @>\s>>  \ca O \\
@V{\b}VV           @|    \\
\ca T'   @>\s'>>  \ca O
\endCD 
$$ 
(and where $\s'$ is the functor defined by obvious analogy with $\s$). 

It is immediate from the commutativity of the diagram that $\b$ is the identity map on objects. Then 
for each $P\in\D$, $\b$ restricts to a group isomorphism 
$$ 
\b_P:Aut_{\ca T}(P)\to Aut_{\ca T'}(P).  
$$ 
In particular, $\b$ maps the identity element $\iota_P$ to $\iota'_P$, where 
$\iota_P=(\1)\e_P$ and where $\iota'_P=(\1)\e'$. Now let $P\leq Q$ in $\D$. Then 
$\iota_{P,Q}$ is the unique (by A.13(d)) $\ca T$-morphism $\g$ having the property that 
$\iota_P\circ\g=\iota_Q$. Similarly, letting $\iota'_{P,Q}$ denote the image of $\1$ under $\e_{P,Q}$, 
then $\iota'_{P,Q}$ is the unique $\ca T'$-morphism $\g'$ such that $\iota'_P\circ\g'=\iota'_Q$.
Thus $(\iota_{P,Q})\b=\iota'_{P,Q}$. Now let  
$\psi:P\to Q$ and $\bar\psi:\bar P\to\bar Q$ be $\ca T$-isomorphisms such that $\psi\up\bar\psi$. 
That it, assume that $\iota_{P,\bar P}\circ\bar\psi=\psi\circ\iota_{Q,\bar Q}$. Applying $\b$ then 
yields $(\psi)\b\up(\bar\psi)\b$, and thus $\b$ induces a mapping $\a:\ca L\to\ca L'$ (on equivalence 
classes of isomorphisms). The product $\Pi(g_1,\cdots,g_n)$ in $\ca L$ is the equivalence class 
of a composite of a sequence of representatives for $(g_1,\cdots,g_n)$, so $\a$ is a homomorphism 
of partial groups. As $\b$ is invertible, $\a$ is then an isomorphism, as required. 
\qed 
\enddemo

\proclaim {Corollary A.26} Let $\ca F$ be a fusion system on the discrete $p$-toral group $S$. Then 
the following conditions are equivalent. 
\roster 

\item "{(1)}" $\ca F$ is order-saturated (as defined in A.5).  

\item "{(2)}" $\ca F$ is saturated (as defined in 8.2), and $Out_{\ca F}(P)$ is finite for every 
subgroup $P\leq S$. 

\item "{(3)}" $\ca F$ is the fusion system of a compact locality. 

\endroster 
\endproclaim

\demo {Proof} If (3) holds then also (1) and (2) hold, by A.7. Thus it now suffices to show that 
(1) implies (3). Assume (1). By [LL] (and [M]) there exists a $p$-local compact group 
$(\ca T,\e,\r)$ over $\ca F$, and A.23 then yields (3). 
\qed 
\enddemo

\enddocument